\documentclass[amstex,12pt, amssymb]{article}

\usepackage{mathtext}
\usepackage[cp1251]{inputenc}
\usepackage[T2A]{fontenc}
\usepackage[dvips]{graphicx}
\usepackage{amsmath}
\usepackage{amssymb}
\usepackage{amsxtra}
\usepackage{latexsym}
\usepackage{ifthen}

\newcounter{lemma}[section]

\newcounter{corollary}[section]

\newcounter{remark}[section]

\newcounter{theorem}[section]

\newcounter{proposition}[section]

\newcounter{example}

\numberwithin{equation}{section}

\begin{document}

\markboth{\centerline{E.~SEVOST'YANOV, S.~SKVORTSOV, N.~ILKEVYCH}}
{\centerline{ON BEHAVIOR ... }}

\def\cc{\setcounter{equation}{0}
\setcounter{figure}{0}\setcounter{table}{0}}

\overfullrule=0pt


\author{{E.~SEVOST'YANOV, S.~SKVORTSOV, N.~ILKEVYCH}\\}

\title{
{\bf ON BEHAVIOR OF A CLASS OF MAPPINGS IN TERMS OF PRIME ENDS}}

\date{\today}
\maketitle

\begin{abstract} The paper
is devoted to the study of mappings with finite distortion, actively
studied recently. For mappings whose inverse satisfy the Poletsky
inequality, the results on boundary behavior in terms of prime ends
are obtained. In particular, it was proved that the families of the
indicated mappings are equicontinuous at the points of the boundary
if a certain function determining the distortion of the module under
the mappings is integrable in a given domain.
\end{abstract}


\section{Introduction} In the theory of quasiconformal mappings,
an important place is occupied by the results on their local and
boundary behavior, see e.g.~\cite[Theorem~19.2]{Va$_1$},
\cite[Theorem~3.17]{MRV}, \cite[Theorem~3.1]{NP$_1$},
\cite[Theorem~3.1]{NP$_2$}, \cite[Theorem~8.9]{Cr} and
\cite[Theorem~3.1, Corollary~3.6]{MRSY$_2$}. Let us mention the
following very important result, see~\cite[Theorem~3.1]{NP$_1$}.

{\bf Theorem (N\"{a}kki--Palka).}{\sl\, Let $\frak{F}$ be a family
of $K$-quasiconformal mappings of a domain $D\ne \overline{{\Bbb
R}^n}$ onto a domain $D^{\,\prime}$ and let either $D$ or
$D^{\,\prime}$ be quasiconformally collared on the boundary. Then
$\frak{F}$ is uniformly equicontinuous if and only if each
$f\in\frak{F}$ can be extended to a continuous mapping of $D$ onto
$D^{\,\prime}$ and $\inf\limits_{\frak{F}} h(f(A))>0$ for some
continuum $A$ in $D.$}

Similar statements can also be obtained for mappings with unbounded
characteristic, which this article is devoted to. For convenience,
in order to separate our studies, we will speak exclusively about
the boundary behavior of maps and do not consider here their
behavior at inner points. We would also like to note that the
publication is devoted to the study of mappings in domains with bad
boundaries. Similar studies have taken place in some of our earlier
papers, see, for example, \cite{SevSkv$_1$} and \cite{Sev}. Unlike
previous articles, the main attention here is paid to the behavior
of homeomorphisms, the inverse of which satisfy Poletsky-type
inequalities. For quasiconformal mappings, consideration of such
mappings does not make sense, since the inverse of quasiconformal
homeomorphisms belong to the same class. The situation changes
drastically if the characteristic of mappings is unbounded. We will
confirm what we said with one of the examples given at the end of
this article.

Recall some definitions (see, for example,~\cite{KR$_1$} and
\cite{KR$_2$}). Let $\omega$ be an open set in ${\Bbb R}^k$,
$k=1,\ldots,n-1$. A continuous mapping
$\sigma\colon\omega\rightarrow{\Bbb R}^n$ is called a {\it
$k$-dimensional surface} in ${\Bbb R}^n$. A {\it surface} is an
arbitrary $(n-1)$-dimensional surface $\sigma$ in ${\Bbb R}^n.$ A
surface $\sigma$ is called {\it a Jordan surface}, if
$\sigma(x)\ne\sigma(y)$ for $x\ne y$. In the following, we will use
$\sigma$ instead of $\sigma(\omega)\subset {\Bbb R}^n,$
$\overline{\sigma}$ instead of $\overline{\sigma(\omega)}$ and
$\partial\sigma$ instead of
$\overline{\sigma(\omega)}\setminus\sigma(\omega).$ A Jordan surface
$\sigma\colon\omega\rightarrow D$ is called a {\it cut} of $D$, if
$\sigma$ separates $D,$ that is $D\setminus \sigma$ has more than
one component, $\partial\sigma\cap D=\varnothing$ and
$\partial\sigma\cap\partial D\ne\varnothing$.

A sequence of cuts $\sigma_1,\sigma_2,\ldots,\sigma_m,\ldots$ in $D$
is called {\it a chain}, if:

(i) $\overline{\sigma_m}\cap\overline{\sigma_{m+1}}=\varnothing$ for
$m\in {\Bbb N};$ (ii) the set $\sigma_{m+1}$ is contained in exactly
one component $d_m$ of the set $D\setminus \sigma_m,$ wherein
$\sigma_{m-1}\subset D\setminus (\sigma_m\cup d_m)$; (iii)
$\bigcap\limits_{m=1}^{\infty}\,d_m=\varnothing.$

According to the definition, a chain of cuts $\{\sigma_m\}$ defines
a chain of domains $d_m\subset D$, such that $\partial\,d_m\cap
D\subset\sigma_m$ and $d_1\supset d_2\supset\ldots\supset
d_m\supset\ldots\, .$ Two chains of cuts  $\{\sigma_m\}$ and
$\{\sigma_k^{\,\prime}\}$ are called {\it equivalent}, if for each
$m=1,2,\ldots$ the domain $d_m$ contains all the domains
$d_k^{\,\prime},$ except for a finite number, and for each
$k=1,2,\ldots$ the domain $d_k^{\,\prime}$ also contains all domains
$d_m,$ except for a finite number.

\begin{figure}
  \centering\includegraphics[width=250pt]{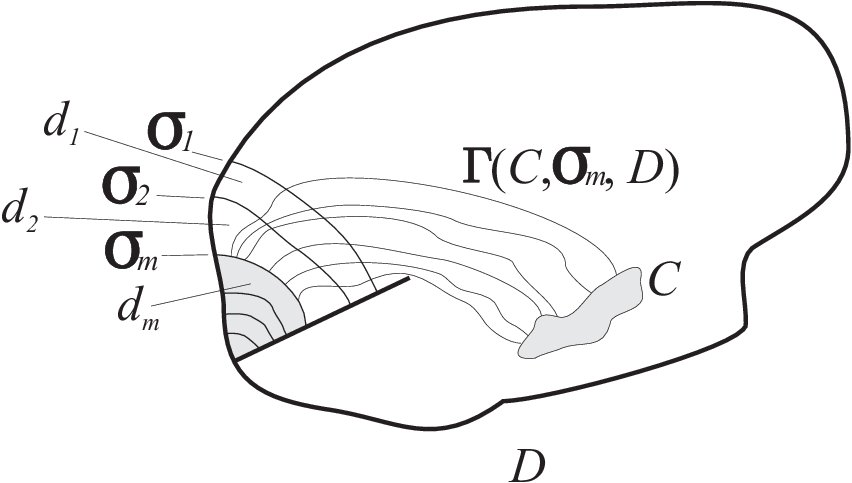}
  \caption{A prime end in some domain}\label{fig1}
 \end{figure}

The {\it end} of the domain $D$ is the class of equivalent chains of
cuts in $D$. Let $K$ be the end of $D$ in ${\Bbb R}^n$, then the set
$I(K)=\bigcap\limits_{m=1}\limits^{\infty}\overline{d_m}$ is called
{\it the impression of the end} $K$. Throughout what follows, as
usual, $\Gamma(E, F, D)$ denotes the family of all paths
$\gamma\colon[a, b]\rightarrow D$ such that $\gamma(a)\in E$ and
$\gamma(b)\in F.$ In what follows, $M$ denotes the modulus of a
family of paths, and the element $dm(x)$ corresponds to a Lebesgue
measure in ${\Bbb R}^n,$ $n\geqslant 2,$ see~\cite{Va$_1$}. For
given sets $E$ and $F$ and a given domain $D$ in $\overline{{\Bbb
R}^n}={\Bbb R}^n\cup \{\infty\},$ we denote by $\Gamma(E, F, D)$ the
family of all paths $\gamma:[0, 1]\rightarrow \overline{{\Bbb R}^n}$
joining $E$ and $F$ in $D,$ that is, $\gamma(0)\in E,$ $\gamma(1)\in
F$ and $\gamma(t)\in D$ for all $t\in (0, 1).$ Following~\cite{Na},
we say that the end $K$ is {\it a prime end}, if $K$ contains a
chain of cuts $\{\sigma_m\}$ such that
$\lim\limits_{m\rightarrow\infty}M(\Gamma(C, \sigma_m, D))=0$ for
some continuum $C$ in $D$ (see Figure~\ref{fig1} for this).
In the following, the following notation is used: the set of prime
ends corresponding to the domain $D,$ is denoted by $E_D,$ and the
completion of the domain $D$ by its prime ends is denoted
$\overline{D}_P.$

Consider the following definition, which goes back to
N\"akki~\cite{Na}, see also~\cite{KR$_1$}--\cite{KR$_2$}. We say
that the boundary of the domain $D$ in ${\Bbb R}^n$ is {\it locally
quasiconformal}, if each point $x_0\in\partial D$ has a neighborhood
$U$ in ${\Bbb R}^n$, which can be mapped by a quasiconformal mapping
$\varphi$ onto the unit ball ${\Bbb B}^n\subset{\Bbb R}^n$ so that
$\varphi(\partial D\cap U)$ is the intersection of ${\Bbb B}^n$ with
the coordinate hyperplane. For a given set $E\subset {\Bbb R}^n,$ we
set $d(E):=\sup\limits_{x, y\in E}|x-y|.$
The sequence of cuts $\sigma_m,$ $m=1,2,\ldots ,$ is called {\it
regular,} if $d(\sigma_{m})\rightarrow 0$ as $m\rightarrow\infty.$
If the end $K$ contains at least one regular chain, then $K$ will be
called {\it regular}. We say that a bounded domain $D$ in ${\Bbb
R}^n$ is {\it regular}, if $D$ can be quasiconformally mapped to a
domain with a locally quasiconformal boundary whose closure is a
compact in ${\Bbb R}^n.$ Note that each prime end of the regular
domain contains a regular chain of cuts, and vice versa, if the
specified property occurs at this end, then it is prime (see
e.g.~\cite[Theorem~5.1]{Na}). We define the closure of a domain with
respect to the space of prime ends by the relation
$\overline{D}_P:=D\cup E_D.$ Note that this space is metric, which
can be demonstrated as follows. If $g:D_0\rightarrow D$ is a
quasiconformal mapping of a domain $D_0$ with a locally
quasiconformal boundary onto some domain $D,$ then for $x, y\in
\overline{D}_P$ we put:
\begin{equation}\label{eq5}
\rho(x, y):=|g^{\,-1}(x)-g^{\,-1}(y)|\,,
\end{equation}
where the element $g^{\,-1}(x),$ $x\in E_D,$ is to be understood as
some (single) boundary point of the domain $D_0.$ The specified
boundary point is well-defined by~\cite[Theorem~4.1]{Na}. It is easy
to verify that~$\rho$ in~(\ref{eq5}) is a metric on
$\overline{D}_P,$ and that the topology on $\overline{D}_P,$ defined
by such a method, does not depend on the choice of the map $g$ with
the indicated property.

\begin{figure}[h]
\centering\includegraphics[width=350pt]{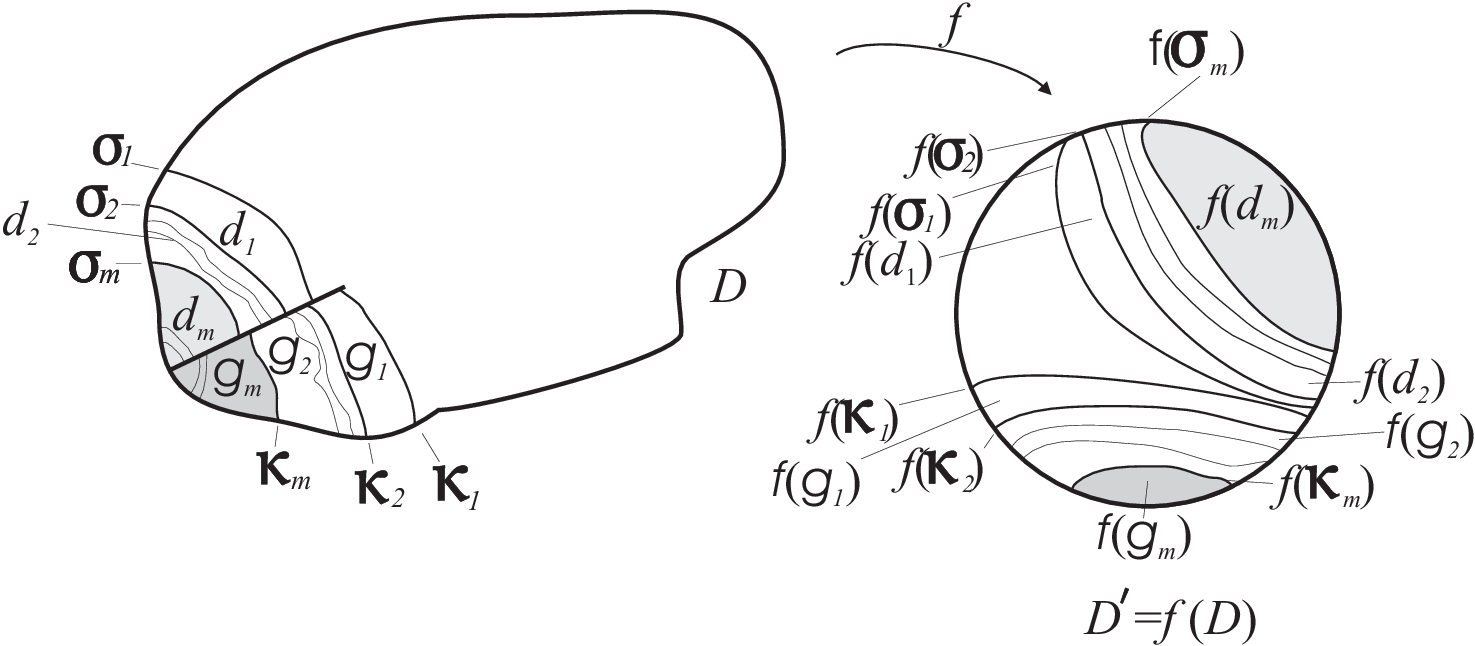} \caption{The
correspondence between the end and the sequence of domains
}\label{fig2}
\end{figure}

We say that a sequence $x_m\in D,$ $m=1,2,\ldots,$ converges to a
prime end of $P\in E_D$ as $m\rightarrow\infty, $ if for any natural
$k\in {\Bbb N}$ all elements of the sequence $x_m$ belong to $d_k$
except for a finite number. Here $d_k$ denotes a sequence of nested
domains corresponding to the definition of the prime end $P.$ Note
that for a homeomorphism of a domain $D$ onto $D^{\,\prime},$ the
end of the domain $D$ uniquely corresponds to some sequence of
nested domains in the image under the mapping. See Figure~\ref{fig2}
for an illustration.
Everywhere below, unless otherwise stated, the boundary and the
closure of a set are understood in the sense of an extended
Euclidean space $\overline{{\Bbb R}^n}.$ Let $x_0\in\overline{D},$
$x_0\ne\infty,$
$$S(x_0,r) = \{
x\,\in\,{\Bbb R}^n : |x-x_0|=r\}\,, S_i=S(x_0, r_i)\,,\quad
i=1,2\,,$$
$$A=A(x_0, r_1, r_2)=\{ x\,\in\,{\Bbb R}^n : r_1<|x-x_0|<r_2\}\,.$$
Let $Q:{\Bbb R}^n\rightarrow {\Bbb R}^n$ be a Lebesgue measurable
function satisfying the condition $Q(x)\equiv 0$ for $x\in{\Bbb
R}^n\setminus D.$ The mapping $f:D\rightarrow \overline{{\Bbb R}^n}$
is called a {\it ring $Q$-mapping at the point $x_0\in
\overline{D}\setminus \{\infty\}$}, if the condition
\begin{equation} \label{eq2*!}
M(f(\Gamma(S_1, S_2, D)))\leqslant \int\limits_{A\cap D} Q(x)\cdot
\eta^n (|x-x_0|)\, dm(x)
\end{equation}
holds for all $0<r_1<r_2<d_0:=\sup\limits_{x\in D}|x-x_0|$ and all
Lebesgue measurable functions $\eta:(r_1, r_2)\rightarrow [0,
\infty]$ such that
\begin{equation}\label{eq8B}
\int\limits_{r_1}^{r_2}\eta(r)dr\geqslant 1\,.
\end{equation}
The mapping of $f$ is called a {\it ring $Q$-mapping in $D,$} if
condition~(\ref{eq2*!}) is satisfied at every point $x_0\in D,$ and
a {\it ring $Q$-mapping in $\overline{D},$} if the
condition~(\ref{eq2*!}) holds at every point $x_0\in\overline{D}.$
With regard to the definition of such mappings, we point to
publications~\cite{RSY} and~\cite{MRSY$_2$}. The class of mappings
satisfying relation (\ref{eq2*!}) contains in itself all conformal
and quasiconformal mappings, as well as many mappings with finite
distortion, see, for example,~\cite[Theorem~2]{Pol},
\cite[Theorem~3.1]{Va$_2$} and~\cite[Theorems~4.6 and
6.10]{MRSY$_1$}. 

Let $(X, d)$ and $\left(X^{{\,\prime}}, d^{\,\prime}\right)$ be
metric spaces with distances $d$ and $d^{\,\prime},$ respectively. A
family $\frak{G}$ of mappings $g:X^{\,\prime}\rightarrow X$  is said
to be {\it equicontinuous at a point} $y_0\in X^{\,\prime},$ if for
every $\varepsilon>0$ there is $\delta=\delta(\varepsilon, y_0)>0$
such that $d(g(y), g(y_0))<\varepsilon$ for all $g\in \frak{G}$ and
$y\in X^{\,\prime}$ with $d^{\,\prime}(y , y_0)<\delta$. The family
$\frak{G}$ is {\it equicontinuous} if $\frak{G}$ is equicontinuous
at every point $y_0\in X^{\,\prime}.$

Everywhere below, unless otherwise stated, $d=\rho$ is one of the
metrics in $\overline{D}_P,$ defined by the relation~(\ref{eq5}),
and $d^{\,\prime}=h$ is a chordal metric defined by formula
\begin{equation}\label{eq1E}
q(x,y)=\frac{|x-y|}{\sqrt{1+{|x|}^2} \sqrt{1+{|y|}^2}}\,,\quad x\ne
\infty\ne y\,, \quad\,q(x,\infty)=\frac{1}{\sqrt{1+{|x|}^2}}\,.
\end{equation}
For a given set $E\subset\overline{{\Bbb R}^n},$ we set
\begin{equation}\label{eq9C}
q(E):=\sup\limits_{x,y\in E}q(x, y)\,.
\end{equation}
The quantity $q(E)$ is called the {\it chordal diameter} of the set
$E.$ The boundary of the domain $D$ is called {\it weakly flat at
the point $x_0,$} if for every number $P>0$ and for every
neighborhood $U$ of this point there is a neighborhood $V$ of point
$x_0$ such that $M(\Gamma(E, F, D))>P$ for arbitrary continua $E$
and $F,$ satisfying conditions $F\cap \partial U\ne\varnothing\ne
F\cap
\partial V.$ The boundary of domain
$D$ is called {\it weakly flat} if it is such at each point of its
boundary.

For a given number $\delta>0,$ domains $D\subset {\Bbb R}^n$ and
$D^{\,\prime}\subset\overline{{\Bbb R}^n},$ $n\geqslant 2,$ a
continuum $A\subset D$ and a Lebesgue measurable function $Q(x):
{\Bbb R}^n\rightarrow [0, \infty]$ such that $Q(x)\equiv 0$ for
$x\not\in D,$ we denote by ${\frak S}_{\delta, A, Q }(D,
D^{\,\prime})$ the family of all homeomorphisms $h$ of
$D^{\,\prime}$ onto $D$ such that the mapping $f=h^{\,-1}$ satisfies
the condition~(\ref{eq2*!}) in $\overline{D},$ while
$q(f(A))\geqslant\delta.$ The following statement is true.

\medskip
\begin{theorem}\label{th2}
{\sl Suppose that $D$ is regular, $D^{\,\prime}$ has a weakly flat
boundary, and any component of $\partial D^{\,\prime}$ is a
non-degenerate continuum. If $Q\in L^1(D),$ then each map $h\in
{\frak S}_{\delta, A, Q }(D, D^{\,\prime})$ extends by continuity to
the map $\overline{h}:\overline{D^{\,\prime}}\rightarrow
\overline{D}_P,$ in addition,
$\overline{h}(D^{\,\prime})=\overline{D}_P,$ and the family ${\frak
S}_{\delta, A, Q }(\overline{D}_P, \overline{D^{\,\prime}}),$
consisting of all extended mappings
$\overline{h}:\overline{D^{\,\prime}}\rightarrow \overline{D}_P,$ is
equicontinuous in $\overline{D^{\,\prime}}.$ }
\end{theorem}

\medskip
\begin{remark}\label{rem1}
The possibility of continuous extension of a homeomorphism
$h:D^{\,\prime}\rightarrow D$ to the mapping
 $\overline{h}:\overline{D^{\,\prime}}\rightarrow
\overline{D}_P$ in Theorem~\ref{th2} may be established similarly
to~\cite[Theorem~6.1]{GRY}; see also~\cite[Theorem~2]{SalSev}. Since
the proof of this result almost literally repeats the reasoning
related to the mentioned publications, we will not give this proof
in the present text.
\end{remark}

\medskip
\section{Preliminaries} Recall that a {\it path} will
be called a continuous mapping $\gamma: I\rightarrow {\Bbb R}^n$ of
a segment, interval or half-interval $I\subset {\Bbb R}$ into
$n$-dimensional Euclidean space ${\Bbb R}^n.$ As usual, the set
$|\gamma|=\{x\in {\Bbb R}^n: \exists\, t\in [a, b]: \gamma(t)=x\}$
is called the {\it locus} of a path $\gamma: I\rightarrow {\Bbb
R}^n.$ We say that the path $\gamma$ lies in the domain $D,$ if its
locus belongs to this domain. We also say that the paths $\gamma_1$
and $\gamma_2$ do not intersect each other if their loci do not
intersect as sets in ${\Bbb R}^n.$ By definition, a prime end $P\in
E_D$ corresponds to a sequence of nested domains $d_m,$ $m\geqslant
1,$ and if $P\in D,$ then we assume that $P$ corresponds to a
sequence of balls $B(P, r_m)$ with radii $r_m\rightarrow 0$ as
$m\rightarrow\infty,$ $r_m>0,$ which lie in the domain of $D$ along
with its closure. Strictly speaking, such a sequence of balls does
not correspond to any prime end in our understanding of the word.
The following statement is true, see
e.g.~\cite[Proposition~1]{SevSkv$_1$}.

\medskip
\begin{proposition}\label{pr1}{\sl\,
Let $n\geqslant 2, $ and let $D$ be a domain in ${\Bbb R}^n$ that is
locally connected on its boundary. Then every two pairs of points
$a\in D, b\in \overline{D}$ and $c\in D, d\in \overline{D}$ can be
joined by non-intersecting paths $\gamma_1:[0, 1]\rightarrow
\overline{D}$ and $\gamma_2:[0, 1]\rightarrow \overline{D}$ so that
$\gamma_i(t)\in D$ for all $t\in (0, 1)$ and all $i=1,2,$ while
$\gamma_1(0)=a,$ $\gamma_1(1)=b,$ $\gamma_2(0)=c$ and
$\gamma_2(1)=d.$}
\end{proposition}

\medskip
The proof of the following statement completely repeats the proof
of~\cite[Theorem~17.10]{Va$_1$}, and therefore is omitted.

\begin{proposition}\label{pr2}
{\sl\, Let $D\subset {\Bbb R}^n$ be a domain with a locally
quasiconformal boundary, then the boundary of this domain is weakly
flat. Moreover, the neighborhood of $U$ in the definition of a
locally quasiconformal boundary can be taken arbitrarily small, and
in this definition $\varphi(x_0)=0.$ }
\end{proposition}

\medskip
The following statement points to the possibility of a
''convenient'' joining of the points of a regular domain by paths.

\begin{lemma}\label{lem1}{\sl\,
Let $D\subset {\Bbb R}^n,$ $n\geqslant 2,$ be a regular domain, and
let $x_m\rightarrow P_1,$ $y_m\rightarrow P_2$ as
$m\rightarrow\infty,$ $P_1, P_2\in \overline{D}_P,$ $P_1\ne P_2.$
Suppose $d_m, g_m,$ $m=1,2,\ldots,$ are sequences of descending
domains, corresponding to $P_1$ and $P_2,$ $d_1\cap
g_1=\varnothing,$ and $x_0, y_0\in D\setminus (d_1\cup g_1).$ Then
there are arbitrarily large $k_0\in {\Bbb N},$ $M_0=M_0(k_0)\in
{\Bbb N}$ and $0<t_1=t_1(k_0), t_2=t_2(k_0)<1$ for which the
following condition is fulfilled: for each $m\geqslant M_0$ there
are
non-intersecting paths
$$\gamma_{1,m}(t)=\quad\left\{
\begin{array}{rr}
\widetilde{\alpha}(t), & t\in [0, t_1],\\
\widetilde{\alpha_m}(t), & t\in [t_1, 1]\end{array}
\right.\,,\quad\gamma_{2,m}(t)=\quad\left\{
\begin{array}{rr}
\widetilde{\beta}(t), & t\in [0, t_2],\\
\widetilde{\beta_m}(t), & t\in [t_2, 1]\end{array}\,, \right.$$
such that:

1) $\gamma_{1, m}(0)=x_0,$ $\gamma_{1, m}(1)=x_m,$ $\gamma_{2,
m}(0)=y_0$ and $\gamma_{2, m}(1)=y_m;$

2) $|\gamma_{1, m}|\cap \overline{g_{k_0}}=\varnothing=|\gamma_{2,
m}|\cap \overline{d_{k_0}};$

3) $\widetilde{\alpha_m}(t)\in d_{k_0}$ for $t\in [t_1, 1]$ and
$\widetilde{\beta_m}(t)\in g_{k_0}$ for $t\in [t_2, 1]$ (see
Figure~\ref{fig3}).
\begin{figure}[h]
\centering\includegraphics[width=200pt]{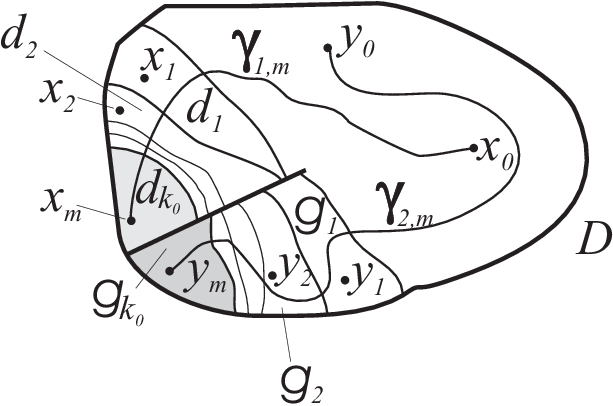} \caption{To
the statement of Lemma~\ref{lem1}}\label{fig3}
\end{figure}
}
\end{lemma}

\begin{proof}
Since, by condition, $D$ is a regular domain, it can be mapped onto
some domain with a locally quasiconformal boundary by (some)
quasiconformal mapping $h:D\rightarrow D_0.$ Note that the domain
$D_0$ is locally connected on its boundary, which follows directly
from the definition of local quasiconformality.

Note that, if $P_1$ and $P_2$ are different prime ends in $D,$ then
$h(P_1)$ and $h(P_2)$ are different prime ends in $D_0.$ Indeed, let
$\sigma_m$ be a sequence of cuts corresponding to the prime
end~$P_1.$ The fact that $h(\sigma_m)$ is also a cut of the domain
$D_0$ is obvious, since $h$ is a homeomorphism. Now we verify that
the sequence $h(\sigma_m),$ $m=1,2,\ldots ,$ is a chain. The
conditions (ii) and (iii) taken from the definition of a chain are
obvious, since $h$ is a homeomorphism. We now verify the condition
(i):
$\overline{h(\sigma_m)}\cap\overline{h(\sigma_{m+1})}=\varnothing$
for $m\in {\Bbb N}.$ Suppose the contrary, namely, that
$\overline{h(\sigma_m)}\cap\overline{h(\sigma_{m+1})}\ne\varnothing$
at least for one $m\in {\Bbb N}.$ Then there is a point $x_0\in
\partial D_0$ such that $x_0\in \overline{h(\sigma_m)}\cap\overline{h(\sigma_{m+1})}.$
Proposition~\ref{pr2} implies that $M(h(\sigma_m), h(\sigma_{m+1}),
D_0)=\infty.$ On the other hand, in view of the definition of the
modulus of families of paths, $M(\Gamma(\sigma_m, \sigma_{m+1},
D))\leqslant l_0^{-n}\cdot m(D)<\infty,$ where $l_0:={\rm
dist}\,(\sigma_m, \sigma_{m+1})>0$ and $m(D)$ is a Lebesgue measure
of $D.$ Here it was also taken into account that the domain $D$ is
bounded, so that $m(D)<\infty.$ Then, due to the quasiconformality
of the mapping $h,$ we have that $M(h(\sigma_m), h(\sigma_{m+1}),
D_0))\leqslant K\cdot M(\sigma_m, \sigma_{m+1}, D_0)<\infty,$ where
$K<\infty$ is some constant. The resulting contradiction refutes the
assumption that
$\overline{h(\sigma_m)}\cap\overline{h(\sigma_{m+1})}\ne\varnothing.$

Thus, the chain of cuts $h(\sigma_m)$ defines some end $h(P_1).$ The
fact that this end is prime also simply follows from the
quasiconformality of the mapping $h.$ Similarly, $h(P_2)$ is a prime
end in $D_0.$

Note that the impressions $I(h(P_1))$ and $I(h(P_2))$ of $h(P_1)$
and $h(P_2)$ are some different points  $a$ and $b$ in $\partial
D_0$ (see~\cite[Theorem~4.1]{Na}). If $P_1$ or $P_2$ are inner
points of $D,$ then $h(P_1)$ (or $h(P_2)$) are inner points of
$D_0,$ which we denote by $a$ or $b,$ respectively. Since by
assumption $x_0, y_0\in D\setminus (d_1\cup g_1),$ then, in
particular, $P_1\ne x_0\ne P_2,$ $P_1\ne y_0\ne P_2.$ This implies
that $a, b, h(x_0)$ and $h(y_0)$ are four different points in
$\overline{D_0},$ at least two of which are inner points of $D_0.$
For the above construction, see Figure~\ref{fig4}.
\begin{figure}[h]
\centering\includegraphics[width=300pt]{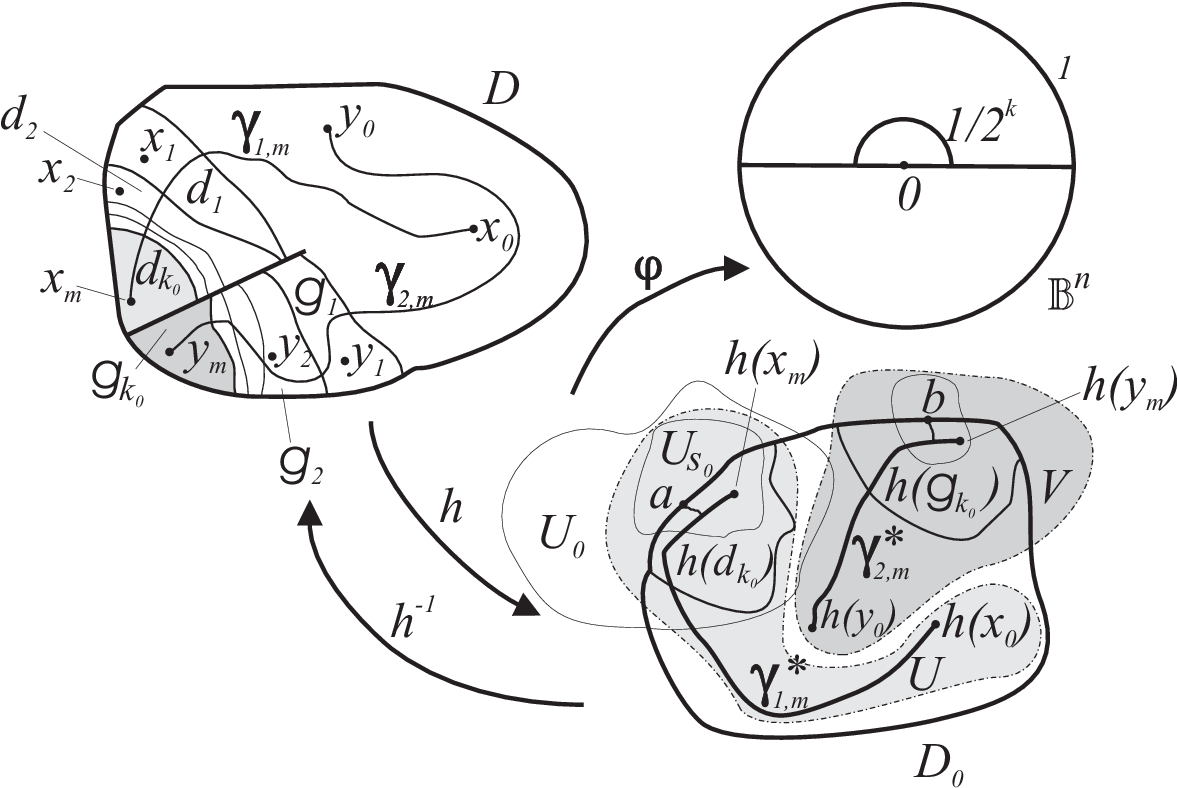} \caption{To
the proof of Lemma~\ref{lem1}}\label{fig4}
\end{figure}
By Proposition~\ref{pr1}, one can join pairs of points $a, h(x_0)$
and $b, h(y_0)$ by disjoint paths $\alpha:[0, 1]\rightarrow
\overline{D_0}$ and $\beta:[0, 1]\rightarrow \overline{D_0}$ so that
$|\alpha|\cap |\beta|=\varnothing,$ $\alpha(t),\beta(t)\in D$ for
all $t\in (0, 1),$ $\alpha(0)=h(x_0),$ $\alpha(1)=a,$
$\beta(0)=h(y_0)$ and $\beta(1)=b.$ Since ${\Bbb R}^n$  is a normal
topological space, the loci $|\alpha|$ and $|\beta|$ have
non-intersecting open neighborhoods $U, V$ such that
\begin{equation}\label{eq3B}
|\alpha|\subset U,\quad |\beta|\subset V\,.
\end{equation}
Here two cases are possible: either $h(P_1)$ is a prime end in
$E_{D_0},$ or a point in  $D_0.$ Let $h(P_1)$ be a prime end in
$E_{D_0}.$ Since $I(h(P_1))=a,$ then there is a number $k_1\in {\Bbb
N}$ such that $\overline{h(d_k)}\subset U$ при $k\geqslant k_1.$ If
$h(P_1)$ is a point of $D,$ then there is also a number $k_1\in
{\Bbb N}$ such that $\overline{h(d_k)}\subset U$ for all $k\geqslant
k_1,$ where $d_k:=B(P_1, r_k),$ $r_k\rightarrow 0$ as
$k\rightarrow\infty$ and $r_k>0.$ In either of these two cases,
$\overline{h(d_k)}\subset U$ for $k\geqslant k_1.$ Similarly, there
is a number $k_2\in {\Bbb N}$ such that $\overline{h(g_k)}\subset V$
for all $k\geqslant k_2.$ Then for $k_0:=\max\{k_1, k_2\}$ we obtain
that
\begin{equation}\label{eq3}
\overline{h(d_k)}\subset U\,,\quad \overline{h(g_k)}\subset V\,,
\quad U\cap V=\varnothing\,, \quad k\geqslant k_0\,.
\end{equation}
Since the sequence $x_m$ converges to $P_1$ as $m\rightarrow\infty,$
then the sequence $h(x_m)$ converges to $a.$ Therefore, there is a
number $m_1\in {\Bbb N}$ such that $h(x_m)\in h(d_{k_0})$ for
$m\geqslant m_1.$ Similarly, since the sequence $y_m$ converges to
$P_2$ as $m\rightarrow\infty,$ then the sequence $h(y_m)$ converges
to $b.$ Therefore, there is a number $m_2\in {\Bbb N}$ such that
$h(y_m)\in h(g_{k_0})$ for $m\geqslant m_2.$ Put $M_0:=\max\{m_1,
m_2\}.$ Show that
\begin{equation}\label{eq2}
|\alpha|\cap h(d_{k_0})\ne\varnothing,\qquad |\beta|\cap
h(g_{k_0})\ne\varnothing\,.
\end{equation}
It suffices to establish the first of these relations, since the
second relation can be proved similarly. If $a=h(P_1)$ is an inner
point of $D_0,$ then this inclusion is obvious. Now suppose that
$h(P_1)$ is a prime end in $E_{D_0}.$ Since the domain $D_0$ has a
locally quasiconformal boundary, there is a sequence of spheres
$S(0,1/2^k),$ $k=0,1,2,\ldots,$ a decreasing sequence of
neighborhoods $U_k$ of the point $a$ and some quasiconformal mapping
$\varphi:U_0\rightarrow {\Bbb B}^n,$ for which $\varphi(U_k)=B(0,
1/2^k),$ $\varphi(\partial U_k\cap D_0)=S(0, 1/2^k)\cap {\Bbb
B}^n_+,$ where ${\Bbb B}^n_+=\{x=(x_1, \ldots, x_n): |x|<1, x_n>0\}$
(see the arguments given in the proof of~\cite[Lemma~3.5]{Na}). Note
that $U_k\cap D_0$ is a domain, since $U_k\cap
D_0=\varphi^{\,-1}(B_+(0, 1/2^k)),$ $B_+(0, 1/2^k)=\{x=(x_1, \ldots,
x_n): |x|<1/2^k, x_n>0\},$ and $\varphi$ is a homeomorphism. In
addition, the sequence of domains $U_k\cap D_0$ corresponds to some
prime end, the impression of which is the point $a,$ and the
corresponding cuts are the sets $\sigma_k:=\partial U_k\cap D_0.$ By
\cite[Theorem~4.1]{Na}, the point $a\in D_0$ corresponds to exactly
one prime end, therefore every domain $h(d_m)$ contains all domains
$U_k\cap D_0,$ except for a finite number, and vice versa. In
particular, there is $s_0\in {\Bbb N}$ such that $U_k\cap D_0\subset
h(d_{k_0})$ for all $k\geqslant s_0.$ Since $a\in |\alpha|,$ there
is $t_1\in (0,1)$ such that $p:=\alpha(t_1)\in U_{s_0}\cap D_0.$ But
then also $p\in h(d_{k_0}),$ since $U_{s_0}\cap D_0\subset
h(d_{k_0}).$ The first relation in~(\ref{eq2}) is proved. As we said
above, the second relation may be proved in exactly the same way.

\medskip
So, let $p:=\alpha(t_1)\in |\alpha|\cap h(d_{k_0}).$ Fix $m\geqslant
M_0$ and join the point $p$ with the point $h(x_m)$ using the path
$\alpha_m:[t_1, 1]\rightarrow h(d_{k_0})$ so that $\alpha_m(t_1)=p,$
$\alpha_m(1)=h(x_m),$ what is possible because $h(d_{k_0})$ is a
domain.
Set
\begin{equation}\label{eq10A}
\gamma^{\,*}_{1,m}(t)=\quad\left\{
\begin{array}{rr}
\alpha(t), & t\in [0, t_1],\\
\alpha_m(t), & t\in [t_1, 1]\end{array} \right.\,.
\end{equation}

Note that the path $\gamma^{\,*}_{1,m}$ completely lies in $U.$
Reasoning similarly, we have the point $t_2\in (0, 1)$ and the point
$q:=\beta(t_2)\in |\beta|\cap h(g_{k_0}).$ Fix $m\geqslant M_0$ and
join the point $q$ with the point $h(y_m)$ using the path
$\beta_m:[t_2, 1]\rightarrow h(g_{k_0})$ so that $\beta_m(t_2)=q,$
$\beta_m(1)=h(y_m),$ that is possible, because $h(g_{k_0})$ is a
domain.
Set
\begin{equation}\label{eq10B}
\gamma^{\,*}_{2,m}(t)=\quad\left\{
\begin{array}{rr}
\beta(t), & t\in [0, t_2],\\
\beta_m(t), & t\in [t_2, 1]\end{array} \right.\,.
\end{equation}
Note that the path $\gamma^{\,*}_{2,m}$ completely lies in $V.$ Set
\begin{equation}\label{eq11A}
\gamma_{1,m}:=h^{\,-1}(\gamma^{\,*}_{1,m})\,,\quad
\gamma_{2,m}:=h^{\,-1}(\gamma^{\,*}_{2,m})\,.
\end{equation}
Note that the paths $\gamma_{1,m}$ and $\gamma_{2,m}$ satisfy all
the conditions of Lemma~\ref{lem1} for $m\geqslant M_0.$ In fact, by
definition, these paths join the points $x_m, x_0$ and $y_m, y_0,$
respectively. The paths $\gamma_{1,m}$ and $\gamma_{2,m}$ do not
intersect, since their images under the mapping $h$ belong to
non-intersecting neighborhoods $U$ and $V,$ respectively.

Note also that $|\gamma_{1,m}|\cap \overline{g_{k_0}}=\varnothing$
for $m\geqslant M_0.$ Indeed, if $x\in |\gamma_{1,m}|\cap
\overline{g_{k_0}},$ then either $x\in |\gamma_{1,m}|\cap g_{k_0}$
or $x\in |\gamma_{1,m}|\cap \partial g_{k_0}.$ In the first case, if
$x\in |\gamma_{1,m}|\cap g_{k_0}$ then $h(x)\in
|\gamma^{\,*}_{1,m}|\cap h(g_{k_0})\subset U\cap h(g_{k_0}),$ which
is impossible due to the relation~(\ref{eq3}). In the second case,
if $x\in |\gamma_{1,m}|\cap \partial g_{k_0},$ then there is a
sequence $z_m\in g_{k_0}$ such that $z_m\rightarrow x$ as
$m\rightarrow\infty.$ Now $h(z_m)\rightarrow h(x)$ as
$m\rightarrow\infty$ and, therefore, $h(x)\in
\overline{h(g_{k_0})}.$ At the same time, $h(x)\in U,$ and this is
impossible by virtue of relation~(\ref{eq3}). Thus, the relation
$|\gamma_{1,m}|\cap \overline{g_{k_0}}=\varnothing$ for $m\geqslant
M_0$ is established.

Similarly, $|\gamma_{2,m}|\cap \overline{d_{k_0}}=\varnothing$ for
$m\geqslant M_0.$ Finally, defining paths $\widetilde{\alpha},$
$\widetilde{\alpha}_m,$ $\widetilde{\beta}$ and
$\widetilde{\beta}_m$ by means of relations
$\widetilde{\alpha}=h^{\,-1}(\alpha),$
$\widetilde{\alpha}_m=h^{\,-1}(\alpha_m),$
$\widetilde{\beta}=h^{\,-1}(\beta)$ and
$\widetilde{\beta}_m=h^{\,-1}(\beta_m),$ we see that these paths
correspond to the construction of $\gamma_{1,m}(t)$ and
$\gamma_{2,m}(t),$ and also satisfy conditions~3) from the
formulation of the lemma. Lemma~\ref{lem1} is proved.~$\Box$
\end{proof}

\medskip
Consider the family of paths joining $\gamma_{1, m}$ and $\gamma_{2,
m}$ from the previous lemma. The following statement contains the
upper estimate of the modulus of the transformed family of paths
under the mapping $f$ with the inequality~(\ref{eq2*!}).

\medskip
\begin{lemma}\label{lem4}
{\sl Let $D\subset {\Bbb R}^n,$ $n\geqslant 2,$ be a regular domain
in ${\Bbb R}^n,$ and let $f:D\rightarrow \overline{{\Bbb R}^n}$ be a
continuous map satisfying the estimate~(\ref{eq2*!}) at every point
$x_0\in\overline{D}$ and some $Q\in L^1(D).$ Then, under the
conditions and notation of Lemma~\ref{lem1}, it is possible to
choose the sequence of domains $d_m$ and the number $k_0$ in such a
way that there exist a constant $0<N=N(k_0, Q, D)<\infty,$
independent of the parameter $m$ and a mapping $f,$ under which
$$M(f(\Gamma(|\gamma_{1, m}|, |\gamma_{2, m}|,
D)))\leqslant N,\qquad m\geqslant M_0=M_0(k_0)\,.$$
}
\end{lemma}
\begin{proof}  By~\cite[Lemma~1]{KR$_2$} the prime
end $P_1$ contains a chain of cuts $\sigma_m$ lying on spheres
$S(\overline{x_0}, r_m)$ such that $\overline{x_0}\in \partial D$
and $r_m\rightarrow 0$ as $m\rightarrow\infty.$ Let $d_m$ be a
sequence of domains corresponding to cuts $\sigma_m.$ Consider
$M_0=M_0(k_0)$ and paths $\gamma_{1,m}$ and $\gamma_{2,m}$
corresponding to this number.

Using the notation of Lemma~\ref{lem1}, we put
$$\varepsilon_0:=\min\{{\rm dist}\,(|\widetilde{\alpha}|, \overline{g_{k_0}}),
{\rm dist}\,(|\widetilde{\alpha}|, |\widetilde{\beta}|)\}>0\,.$$
Now, consider covering of $|\widetilde{\alpha}|$ of the following
type: $\bigcup\limits_{x\in |\widetilde{\alpha}|}B(x,
\varepsilon_0/4).$ Since $|\widetilde{\alpha}|$ is compact in $D,$
there are $i_1,\ldots, i_{N_0}$ such that
$|\widetilde{\alpha}|\subset \bigcup\limits_{i=1}^{N_0} B(z_i,
\varepsilon_0/4),$ where $z_i\in |\widetilde{\alpha}|$ for
$1\leqslant i\leqslant N_0.$ Taking into
account~\cite[Theorem~1.I.5.46]{Ku}, it is easy to verify that
$$\Gamma(|\widetilde{\alpha}|, |\gamma_{2, m}|, D)\subset\bigcup\limits_{i=1}^{N_0}
\Gamma(S(z_i, \varepsilon_0/4), S(z_i, \varepsilon_0/2), D)\,.$$

Putting
$$\eta(t)= \left\{
\begin{array}{rr}
4/\varepsilon_0, & t\in [\varepsilon_0/4, \varepsilon_0/2],\\
0,  &  t\not\in [\varepsilon_0/4, \varepsilon_0/2]
\end{array}
\right. \,,$$
we observe that the function~$\eta$ satisfies relation~(\ref{eq8B}).
Then, by the definition of a ring $Q$-map in~(\ref{eq2*!}) and
taking into account the semi-additivity of the modulus of families
of paths, see~\cite[Theorem~6.2]{Va$_1$}, we obtain that
$$M(f(\Gamma(|\widetilde{\alpha}|, |\gamma_{2, m}|, D)))\leqslant$$
\begin{equation}\label{eq1}
\leqslant \sum\limits_{i=1}^{N_0} M(f(\Gamma(S(z_i,
\varepsilon_0/4), S(z_i, \varepsilon_0/2), D)))\leqslant
\frac{N_04^n\Vert Q\Vert_1}{\varepsilon^n_0}\,, m\geqslant M_0\,,
\end{equation}
where $\Vert Q\Vert_1=\int\limits_DQ(x)\,dm(x).$
On the other hand, by~\cite[Theorem~1.I.5.46]{Ku}, we obtain that
$$\Gamma(|\widetilde{\alpha_m}|, |\gamma_{2, m}|, D)\leqslant
\Gamma(S(\overline{x_0}, r_{k_1}), S(\overline{x_0}, r_{k_2}),
D)\,.$$
Arguing as above, choosing an admissible function
$$\eta(t)= \left\{
\begin{array}{rr}
1/(r_{k_2}-r_{k_1}), & t\in [r_{k_1}, r_{k_2}],\\
0,  &  t\not\in [r_{k_1}, r_{k_2}]
\end{array}
\right. \,,$$
we obtain that
$$M(f(\Gamma(|\widetilde{\alpha_m}|, |\gamma_{2, m}|, D)))\leqslant$$
\begin{equation}\label{eq4D}
\leqslant M(f(\Gamma(S(\overline{x_0}, r_{k_1}), S(\overline{x_0},
r_{k_2}), D)))\leqslant \frac{\Vert
Q\Vert_1}{(r_{k_2}-r_{k_1})^n}\,, m\geqslant M_0\,.
\end{equation}
Now note that $$\Gamma(|\gamma_{1, m}|, |\gamma_{2, m}|, D)\subset
\Gamma(|\widetilde{\alpha_m}|, |\gamma_{2, m}|, D)\cup
\Gamma(|\widetilde{\alpha}|, |\gamma_{2, m}|, D)\,.$$ In this case,
from~(\ref{eq1}) and~(\ref{eq4D}), taking into account the
semi-additivity of the modulus of families of paths, we obtain:
$$M(f(\Gamma(|\gamma_{1, m}|, |\gamma_{2, m}|,
D)))\leqslant
\left(\frac{N_04^n}{\varepsilon^n_0}+\frac{1}{(r_{k_2}-r_{k_1})^n}\right)\Vert
Q\Vert_1\,,\quad m\geqslant M_0\,.$$
The right side of the last relation does not depend on $m,$ so that
we can put
$N:=\left(\frac{N_04^n}{\varepsilon^n_0}+\frac{1}{(r_{k_2}-r_{k_1})^n}\right)\Vert
Q\Vert_1.$ Lemma~\ref{lem4} is completely proved.~$\Box$
\end{proof}

\medskip
The following statement indicates that for some wide class of
mappings fixing the diameter of the image of a certain
non-degenerate continuum, the image of this continuum cannot be
close to the boundary of the corresponding domain under these
mappings.  Note that similar statements were previously known for
quasiconformal mappings, see, for example, \cite[Theorems~21.13 and
21.14]{Va$_1$}. We may also point to our recent result on this,
see~\cite[Lemma~2(v)]{SevSkv$_1$}.

\medskip
\begin{lemma}\label{lem3}{\sl\,
Let $n\geqslant 2,$ let $D$ be a regular domain in ${\Bbb R}^n,$ and
let $D^{\,\prime}$ be some domain in $\overline{{\Bbb R}^n}.$
Suppose that $D^{\,\prime}$ has a weakly flat boundary, $Q\in
L^1(D)$ and, moreover, no connected component of the set $\partial
D^{\,\prime}$ does not degenerate into a point. Let
$f_m:D\rightarrow D^{\,\prime}$ be a sequence of homeomorphisms of
$D$ onto $D^{\,\prime},$ satisfying the relation~(\ref{eq2*!}) in
$D$ with the same function $Q.$

Suppose also that there is a continuum $A\subset D$ and a number
$\delta>0$ such that $q(f_m(A))\geqslant \delta>0$ for all
$m=1,2,\ldots ,$ where, as usual, $q(f_m(A))$ is defined
by~(\ref{eq9C}). Then there is $\delta_1>0$ such that
$$q(f_m(A),
\partial D^{\,\prime})>\delta_1>0\quad \forall\,\, m\in {\Bbb
N}\,,$$
where $q(f_m(A),
\partial D^{\,\prime})=\inf\limits_{x\in f_m(A), y\in \partial D^{\,\prime}}q(x, y).$}
\end{lemma}

\medskip
\begin{proof}
We carry out the proof by contradiction. Suppose that the conclusion
of the lemma is not true. Then for each $k\in {\Bbb N}$ there is
some number $m=m_k$ such that $q(f_{m_k}(A),
\partial D^{\,\prime})<1/k.$
Of course, we can assume that the sequence $m_k$ increases on $k.$
Since $\overline{{\Bbb R}^n}$ is compact, the set $\partial
D^{\,\prime}$ is also compact in extended Euclidean space. Note that
the set $f_{m_k}(A)$ is compact as a continuous image of a compact
set $A\subset D$ under the mapping~$f_{m_k}.$ In this case, there
are elements $x_k\in f_{m_k}(A)$ and $y_k\in
\partial D^{\,\prime}$ such that $q(f_{m_k}(A),
\partial D^{\,\prime})=q(x_k, y_k)<1/k$ (see~Figure~\ref{fig5}).
\begin{figure}
  \centering\includegraphics[width=300pt]{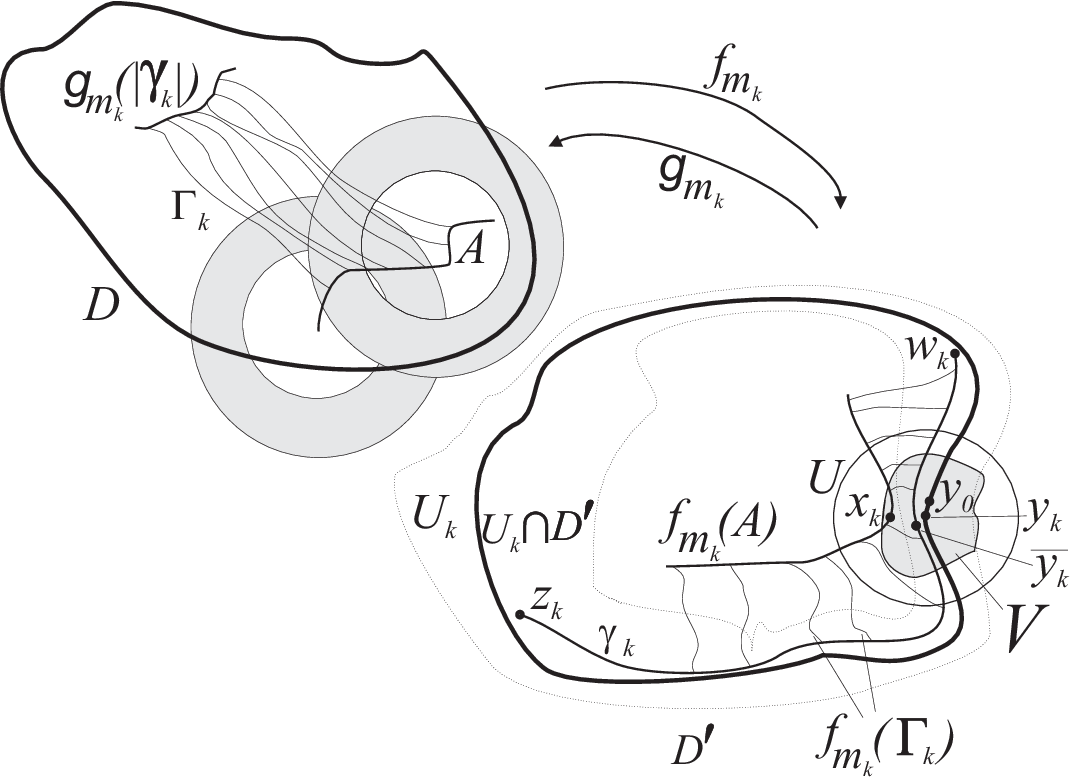}
  \caption{To
the proof of Lemma~\ref{lem3}.}\label{fig5}
 \end{figure}
Since $\partial D^{\,\prime}$ is a compact set, we may assume that
$y_k\rightarrow y_0\in \partial D^{\,\prime}$ as $k\rightarrow
\infty;$ then also
%
$$x_k\rightarrow y_0\in \partial D^{\,\prime},\quad k\rightarrow
\infty\,.$$
%
Let $K_0$ be a connected component of the set $\partial
D^{\,\prime},$ containing $y_0.$ Obviously, $K_0$ is a nondegenerate
continuum in $\overline{{\Bbb R}^n}.$ Since $D^{\,\prime}$ has a
weakly flat boundary, the mapping $g_{m_k}:=f_{m_k}^{\,-1}$ can be
extended to a continuous mapping
$\overline{g}_{m_k}:\overline{D^{\,\prime}}\rightarrow \overline{D}$
(see~Remark~\ref{rem1}). Moreover, $\overline{g}_{m_k}$ is uniformly
continuous on the set $\overline{D^{\,\prime}}$ for every fixed $k,$
because the mapping $\overline{g}_{m_k}$ is continuous on the
compact set $\overline{D^{\,\prime}}.$ Let $\rho$ be one of the
metrics in $E_D=\overline{D}_P\setminus D,$ defined in~(\ref{eq5}),
and let $g:D_0\rightarrow D$ be a quasiconformal mapping of some
domain $D_0$ with locally quasiconformal boundary corresponding to
the definition of the metric $\rho$ in~(\ref{eq5}). In this case,
for each $\varepsilon>0$ there is
$\delta_k=\delta_k(\varepsilon)<1/k$ such that
\begin{equation}\label{eq3A}
\rho(\overline{g}_{m_k}(x),\overline{g}_{m_k}(x_0))<\varepsilon
\quad \forall\,\, x,x_0\in \overline{D^{\,\prime}},\quad q(x,
x_0)<\delta_k\,, \quad \delta_k<1/k\,.
\end{equation}
Choose $\varepsilon>0$ such that
\begin{equation}\label{eq5A}
\varepsilon<(1/2)\cdot {\rm  dist}\,(\partial D_0, g^{\,-1}(A))\,,
\end{equation}
where $A$ is a continuum from the conditions of the lemma. Denote
$B_q(x_0, r)=\{x\in \overline{{\Bbb R}^n}: q(x, x_0)<r\}.$ For a
given $k\in {\Bbb N},$ we set
$$B_k:=\bigcup\limits_{x_0\in K_0}B_q(x_0, \delta_k)\,,\quad k\in {\Bbb N}\,.$$
Since the set $B_k$ is a neighborhood of the continuum $K_0,$ due
to~\cite[Lemma~2.2]{HK} there is a neighborhood $U_k$ of the set
$K_0,$ such that $U_k\subset B_k$ and $U_k\cap D^{\,\prime}$ is
connected. Without loss of generality, we may assume that $U_k$ is
an open set, so $U_k\cap D^{\,\prime}$ is also path connected
(see~\cite[Proposition~13.1]{MRSY$_3$}). Let $q(K_0)=m_0,$ where the
chordal diameter $q(K_0)$ of the set $K_0$ is defined by the
relation~(\ref{eq9C}). In this case, there are $z_0, w_0\in K_0$
such that $q(K_0)=q(z_0, w_0)=m_0.$ So, there are sequences
$\overline{y_k}\in U_k\cap D^{\,\prime},$ $z_k\in U_k\cap
D^{\,\prime}$ and $w_k\in U_k\cap D^{\,\prime}$ such that
$z_k\rightarrow z_0,$ $\overline{y_k}\rightarrow y_0$ and
$w_k\rightarrow w_0$ as $k\rightarrow\infty.$ We may assume that
\begin{equation}\label{eq2A}
q(z_k, w_k)>m_0/2\quad \forall\,\, k\in {\Bbb N}\,.
\end{equation}
Since the set $U_k\cap D^{\,\prime}$ is path-connected, we can
sequentially join the points $z_k,$ $\overline{y_k}$ and $w_k$ using
some path $\gamma_k\in U_k\cap D^{\,\prime}.$ As usual, we denote by
$|\gamma_k|$ the locus of the path $\gamma_k$ in the domain
$D^{\,\prime}.$ Then $g_{m_k}(|\gamma_k|)$ is a compact set in the
domain $D.$ If $x\in|\gamma_k|,$ then there is $x_0\in K_0$ such
that $x\in B(x_0, \delta_k).$ Put $\omega\in A\subset D.$ Since
$x\in|\gamma_k|$ and, moreover, $x$ is an inner point of the domain
$D^{\,\prime},$ we can write here $g_{m_k}(x)$ instead of
$\overline{g}_{m_k}(x).$ By the relations~(\ref{eq3A})
and~(\ref{eq5A}), as well as by the triangle inequality, we obtain
that for sufficiently large $k\in {\Bbb N},$
$$\rho(g_{m_k}(x), \omega)\geqslant
\rho(\omega, \overline{g}_{m_k}(x_0))- \rho(\overline{g}_{m_k}(x_0),
g_{m_k}(x))\geqslant$$
\begin{equation}\label{eq4}
\geqslant {\rm dist}\,(\partial D_0, g^{\,-1}(A))-(1/2)\cdot {\rm
dist}\,(\partial D_0, g^{\,-1}(A))=(1/2)\cdot {\rm dist}\,(\partial
D_0, g^{\,-1}(A))>\varepsilon\,,
\end{equation}
where ${\rm dist}(\partial D_0, g^{\,-1}(A)):=\inf\limits_{x\in
\partial D_0, y\in g^{\,-1}(A)} |x-y|.$ Taking $\inf$ in~(\ref{eq4})
over all $x\in |\gamma_k|$ and $\omega\in A,$ we obtain that
\begin{equation}\label{eq5B}
\rho(g_{m_k}(|\gamma_k|), A):=\inf\limits_{x\in g_{m_k}(|\gamma_k|),
y\in A}\rho(x, y)>\varepsilon, \quad\forall\,\, k=1,2,\ldots \,.
\end{equation}
We now show that there exists $\varepsilon_1>0$ such that
\begin{equation}\label{eq6}
{\rm dist}\,(g_{m_k}(|\gamma_k|), A)>\varepsilon_1, \quad\forall\,\,
k=1,2,\ldots \,,
\end{equation}
where ${\rm dist},$ as usual, denotes the Euclidean distance between
the sets $A, B\subset{\Bbb R}^n.$ Indeed, let~(\ref{eq6}) be
violated, then for the number $\varepsilon_l=1/l,$ $l=1,2,\ldots$
there are $\xi_l\in |\gamma_{k_l}|$ and $\zeta_l\in A$ such that
\begin{equation}\label{eq7A}
|g_{m_{k_l}}(\xi_l)-\zeta_l|<1/l\,,\quad l=1,2,\ldots \,.
\end{equation}
Without loss of generality, we may assume that the sequence $k_l,$
$l=1,2,\ldots,$ is increasing. Since~$A$ is compact, we may assume
that the sequence $\zeta_l$ converges to $\zeta_0\in A$ as
$l\rightarrow\infty.$  By the triangle inequality and
from~(\ref{eq7A}) it follows that
\begin{equation}\label{eq8A}
|g_{m_{k_l}}(\xi_l)-\zeta_0|\rightarrow 0\,,\quad
l\rightarrow\infty\,.
\end{equation}
On the other hand, we recall that $\rho(g_{m_k}(x),
\omega)=|g^{\,-1}(g_{m_k}(x))-g^{\,-1}(\omega)|,$ where
$g:D_0\rightarrow D$ is some quasiconformal mapping of $D_0$ onto
$D,$ see~(\ref{eq5}). In particular, $g^{\,-1}$ is continuous in
$D,$ therefore, by the triangle inequality and~(\ref{eq8A}), we
obtain that
$$|g^{\,-1}(g_{m_{k_l}}(\xi_l))-g^{\,-1}(\zeta_l)|\leqslant$$
\begin{equation}\label{eq9A}
\leqslant
|g^{\,-1}(g_{m_{k_l}}(\xi_l))-g^{\,-1}(\zeta_0)|+|g^{\,-1}(\zeta_0)-g^{\,-1}(\zeta_l)|\rightarrow
0,\quad l\rightarrow\infty\,.\end{equation}
However, by definition $\rho$ and from~(\ref{eq9A}) it follows that
$$\rho(g_{m_{k_l}}(|\gamma_{k_l}|), A)\leqslant \rho(g_{m_{k_l}}(\xi_l), \zeta_l)=
|g^{\,-1}(g_{m_{k_l}}(\xi_l))-g^{\,-1}(\zeta_l)|\rightarrow 0, \quad
l\rightarrow\infty\,,$$
which contradicts~(\ref{eq5B}). The resulting contradiction
indicates the validity of~(\ref{eq6}).

\medskip
We cover the continuum $A$ with the help of balls $B(x,
\varepsilon_1/4),$ $x\in A.$ Since $A$ is a compact set, we may
assume that $A\subset \bigcup\limits_{i=1}^{M_0}B(x_i,
\varepsilon_1/4),$ $x_i\in A,$ $i=1,2,\ldots, M_0,$ $1\leqslant
M_0<\infty.$ By definition, $M_0$ depends only on $A,$ in
particular, $M_0$ does non depend on $k.$ We set
\begin{equation}\label{eq5C}
\Gamma_k:=\Gamma(A, g_{m_k}(|\gamma_k|), D)\,.
\end{equation}
Note that
\begin{equation}\label{eq6C}
\Gamma_k=\bigcup\limits_{i=1}^{M_0}\Gamma_{ki}\,,
\end{equation}
where $\Gamma_{ki}$ consists of all paths $\gamma:[0, 1]\rightarrow
D,$ belonging to the family $\Gamma_k,$ such that $\gamma(0)\in
B(x_i, \varepsilon_1/4)$ and $\gamma(1)\in g_{m_k}(|\gamma_k|).$ We
now show that
\begin{equation}\label{eq7C}
\Gamma_{ki}>\Gamma(S(x_i, \varepsilon_1/4), S(x_i, \varepsilon_1/2),
A(x_i, \varepsilon_1/4, \varepsilon_1/2))\,.
\end{equation}
Indeed, let $\gamma\in \Gamma_{ki},$ in other words, $\gamma:[0,
1]\rightarrow D,$ $\gamma(0)\in B(x_i, \varepsilon_1/4)$ and
$\gamma(1)\in g_{m_k}(|\gamma_k|).$ By~(\ref{eq6}), $|\gamma|\cap
B(x_i, \varepsilon_1/4)\ne\varnothing\ne |\gamma|\cap (D\setminus
B(x_i, \varepsilon_1/4)).$ Therefore, by~\cite[Theorem~1.I.5.46]{Ku}
there is $0<t_1<1$ with the condition $\gamma(t_1)\in S(x_i,
\varepsilon_1/4).$ We can assume that $\gamma(t)\not\in B(x_i,
\varepsilon_1/4)$ for $t>t_1.$ Put $\gamma_1:=\gamma|_{[t_1, 1]}.$
By~(\ref{eq6}), $|\gamma_1|\cap B(x_i,
\varepsilon_1/2)\ne\varnothing\ne |\gamma_1|\cap (D\setminus B(x_i,
\varepsilon_1/2)).$ Thus, by~\cite[Theorem~1.I.5.46]{Ku} there is
$t_1<t_2<1$ with $\gamma(t_2)\in S(x_i, \varepsilon_1/2).$ We may
assume that $\gamma(t)\in B(x_i, \varepsilon_1/2)$ for $t<t_2.$ Put
$\gamma_2:=\gamma|_{[t_1, t_2]}.$ Then, the path $\gamma_2$ is a
subpath of $\gamma,$ which belongs to the family $\Gamma(S(x_i,
\varepsilon_1/4), S(x_i, \varepsilon_1/2), A(x_i, \varepsilon_1/4,
\varepsilon_1/2)).$ Thus, the relation~(\ref{eq7C}) is established.

\medskip
Further reasoning is based, as before, on the successful choice of
an admissible function $\eta.$ Put
$$\eta(t)= \left\{
\begin{array}{rr}
4/\varepsilon_1, & t\in [\varepsilon_1/4, \varepsilon_1/2],\\
0,  &  t\not\in [\varepsilon_1/4, \varepsilon_1/2]\,.
\end{array}
\right. $$
Note that $\eta$ satisfies~(\ref{eq8B}) for $r_1=\varepsilon_1/4$
and $r_2=\varepsilon_1/2.$ Then, according to the definition of a
ring $Q$-homeomorphism at $x_i,$ we obtain that
\begin{equation}\label{eq8C}
M(f_{m_k}(\Gamma(S(x_i, \varepsilon_1/4), S(x_i, \varepsilon_1/2)),
A(x_i, \varepsilon_1/4, \varepsilon_1/2)))\leqslant
(4/\varepsilon_1)^n\cdot\Vert Q\Vert_1<c<\infty\,,
\end{equation}
where $c$ is some positive constant and $\Vert Q\Vert_1$ is
$L_1$-norm of the function $Q$ in $D.$ By~(\ref{eq6C}), (\ref{eq7C})
and (\ref{eq8C}), using the subadditivity of modulus, we obtain that
\begin{equation}\label{eq4B}
M(f_{m_k}(\Gamma_k))\leqslant
\frac{4^nM_0}{\varepsilon_1^n}\int\limits_DQ(x)\,dm(x)\leqslant
c\cdot M_0<\infty\,.
\end{equation}
Let us show that the estimate~(\ref{eq4B}) contradicts the condition
of the weak flatness of the boundary of the domain $D^{\,\prime}.$
Let $U:=B_q(y_0, r_0)=\{y\in \overline{{\Bbb R}^n}: q(y,
y_0)<r_0\},$ where $0<r_0<\min\{\delta/4, m_0/4\},$ $\delta$ is the
number from the condition of the lemma and $q(K_0)=m_0.$ (Here, as
usual, $q(K_0)$ denotes the chordal diameter of the set $E=K_0,$
defined by the formula~(\ref{eq9C})). Note that $|\gamma_k|\cap
U\ne\varnothing\ne |\gamma_k|\cap (D^{\,\prime}\setminus U)$ for
sufficiently large $k\in{\Bbb N},$ since $q(|\gamma_k|)>
m_0/2>m_0/4$ by~(\ref{eq2A}), in addition, $\overline{y_k}\in
|\gamma_k|$ and $\overline{y_k}\rightarrow y_0$ as
$k\rightarrow\infty.$ Similarly, $f_{m_k}(A)\cap U\ne\varnothing\ne
f_{m_k}(A)\cap (D^{\,\prime}\setminus U).$ Since $|\gamma_k|$ and
$f_{m_k}(A)$ are continua, we obtain that
\begin{equation}\label{eq8}
f_{m_k}(A)\cap \partial U\ne\varnothing, \quad|\gamma_k|\cap
\partial U\ne\varnothing\,,
\end{equation}
see~\cite[Theorem~1.I.5.46]{Ku}. For a given $P>0,$ let $V\subset U$
be a neighborhood of the point $y_0,$ corresponding to the
definition of a weakly flat boundary. Then we have that
\begin{equation}\label{eq9}
M(\Gamma(E, F, D^{\,\prime}))>P
\end{equation}
for any continua $E, F\subset D^{\,\prime}$ with $E\cap
\partial U\ne\varnothing\ne E\cap \partial V$ and $F\cap \partial
U\ne\varnothing\ne F\cap \partial V.$
Observe that
\begin{equation}\label{eq10}
f_{m_k}(A)\cap \partial V\ne\varnothing, \quad|\gamma_k|\cap
\partial V\ne\varnothing
\end{equation}
for sufficiently large $k\in {\Bbb N}.$ Indeed, $\overline{y_k}\in
|\gamma_k|,$ $x_k\in f_{m_k}(A),$ where $x_k,
\overline{y_k}\rightarrow y_0\in V$ as $k\rightarrow\infty.$
Therefore, $|\gamma_k|\cap V\ne\varnothing\ne f_{m_k}(A)\cap V$ for
large $k\in {\Bbb N}.$ In addition, we have that $q(V)\leqslant
q(U)\leqslant 2r_0<m_0/2.$ By~(\ref{eq2}), $q(|\gamma_k|)>m_0/2,$
therefore, $|\gamma_k|\cap (D^{\,\prime}\setminus V)\ne\varnothing.$
Thus, by~\cite[Theorem~1.I.5.46]{Ku}, $|\gamma_k|\cap\partial
V\ne\varnothing.$ Similarly, $q(V)\leqslant q(U)\leqslant
2r_0<\delta/2.$ Since $q(f_{m_k}(A))>\delta,$ we obtain that
$f_{m_k}(A)\cap (D^{\,\prime}\setminus V)\ne\varnothing.$
By~\cite[Theorem~1.I.5.46]{Ku}, we have that $f_{m_k}(A)\cap
\partial V\ne\varnothing.$ Thus, the relation~(\ref{eq10}) is established.

\medskip
By~(\ref{eq8}), (\ref{eq9}) and (\ref{eq10}), we obtain that
\begin{equation}\label{eq11}
M(\Gamma(f_{m_k}(A), |\gamma_k|, D^{\,\prime}))>P\,.
\end{equation}
Note that $\Gamma(f_{m_k}(A), |\gamma_k|,
D^{\,\prime})=f_{m_k}(\Gamma(A, g_{m_k}(|\gamma_k|),
D))=f_{m_k}(\Gamma_k).$ Therefore, the relation~(\ref{eq11}) can be
written as
$$M(\Gamma(f_{m_k}(A), g_{m_k}(|\gamma_k|), D))=M(f_{m_k}(\Gamma_k))>P\,.$$
The relation obtained above contradicts the estimate~(\ref{eq4B}).
The resulting contradiction means that the above assumption
$q(f_{m_k}(A),
\partial D^{\,\prime})<1/k$ was
incorrect. The proof of the lemma is complete.~~$\Box$
\end{proof}

\medskip
\section{\bf Proof of Theorem~\ref{th2}}
\setcounter{equation}{0}

For the continuous extension of the mapping~$h\in {\frak S}_{\delta,
A, Q }(D, D^{\,\prime})$ to the boundary of the domain
$D^{\,\prime},$ see Remark~\ref{rem1}. The equicontinuity of~${\frak
S}_{\delta, A, Q }(D, D^{\,\prime})$ is the result
of~\cite[Theorem~1.1]{SevSkv$_2$}.

\medskip
We show the equicontinuity ${\frak S}_{\delta, A, Q
}(\overline{D}_P, \overline{D^{\,\prime}})$ on
$E_{D^{\,\prime}}=\overline{D}_P\setminus D.$

\medskip
We carry out the proof by contradiction. Suppose there are a point
$z_0\in \partial D^{\,\prime},$ the number $\varepsilon_0>0,$ the
sequences $z_m\in \overline{D^{\,\prime}},$ $z_m\rightarrow z_0$ as
$m\rightarrow\infty$ and $\overline{h}_m\in {\frak S}_{\delta, A, Q
}(\overline{D}_P, \overline{D^{\,\prime}})$ such that
\begin{equation}\label{eq12}
\rho(\overline{h}_m(z_m),
\overline{h}_m(z_0))\geqslant\varepsilon_0,\quad m=1,2,\ldots ,
\end{equation}
where $\rho$ is one of the metrics in $\overline{D}_P,$ defined by
the formula~(\ref{eq5}). Since $h_m=\overline{h}_m|_{D^{\,\prime}}$
extends by continuity to the boundary of $D^{\,\prime},$ we may
assume that $z_m\in D$ and, in addition, there is another sequence
$z^{\,\prime}_m\in \overline{D^{\,\prime}},$
$z^{\,\prime}_m\rightarrow z_0$ as $m\rightarrow\infty,$ such that
$\rho(h_m(z_m), \overline{h}_m(z_0))\rightarrow 0$ as
$m\rightarrow\infty.$ Then from~(\ref{eq12}) it follows that
\begin{equation}\label{eq13}
\rho(h_m(z_m), h_m(z^{\,\prime}_m))\geqslant\varepsilon_0/2,\quad
m\geqslant m_0\,.
\end{equation}
Since the domain $D$ is regular, the space $\overline{D}_P$ is
compact. Therefore, we may assume that the sequences $h_m(z_m)$ и
$\overline{h}_m(z_0)$ converge as $m\rightarrow\infty$ to some
elements $P_1, P_2\in \overline{D}_P,$ $P_1\ne P_2.$ Let $d_m$ and
$g_m$ be sequences of descending domains corresponding to prime ends
$P_1$ and $P_2,$ respectively. By~\cite[Lemma~1]{KR$_2$} we may
consider that the cuts $\sigma_m$ corresponding to domains $d_m,$
$m=1,2,\ldots, $ belong to spheres $S(\overline{x_0}, r_m)$ so that
$\overline{x_0}\in
\partial D$ and $r_m\rightarrow 0$ as $m\rightarrow\infty.$
Choose $x_0, y_0\in A$ so that $x_0\ne y_0$ and $x_0\ne P_1\ne y_0,$
where the continuum $A\subset D$ is from conditions of Theorem~
\ref{th2}. Without loss of generality, we may assume that $d_1\cap
g_1=\varnothing$ and $x_0, y_0\not\in d_1\cup g_1.$

\begin{figure}
  \centering\includegraphics[width=450pt]{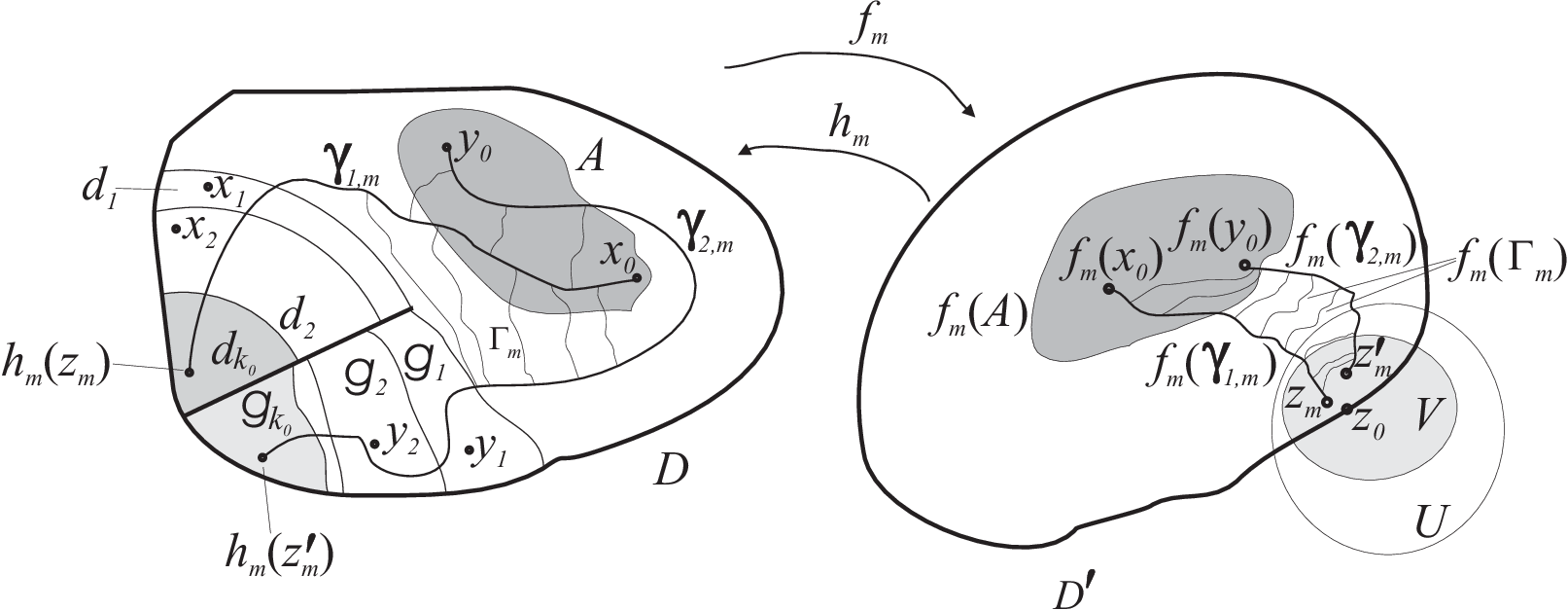}
  \caption{To the proof of Theorem~\ref{th2}.}\label{fig6}
 \end{figure}

\medskip
By Lemmas~\ref{lem1} and~\ref{lem4}, there are disjoint paths
$\gamma_{1,m}:[0, 1]\rightarrow D$ and $\gamma_{2,m}:[0,
1]\rightarrow D,$ the number $M_0>0$ and the number $N> 0$ such that
$\gamma_{1, m}(0)=x_0,$ $\gamma_{1, m}(1)=h_m(z_m),$ $\gamma_{2,
m}(0)=y_0,$ $\gamma_{2, m}(0)=h_m(z^{\,\prime}_m),$ wherein,
\begin{equation}\label{eq15}
M(f_m(\Gamma_m))\leqslant N\,, m\geqslant M_0\,,
\end{equation}
where $f_m:=h^{\,-1}_m,$ $\Gamma_m:=\Gamma(|\gamma_{1, m}|,
|\gamma_{2, m}|, D).$ See Figure~\ref{fig6} for further explanation
of the construction of the proof.
On the other hand, by Lemma~\ref{lem3} there is a number
$\delta_1>0$ such that $q(f_{m}(A), \partial
D^{\,\prime})>\delta_1>0,$ $m=1,2,\ldots \,.$ From this we obtain
that
$$q(f_m(|\gamma_{1, m}|))\geqslant q(z_m, f_m(x_0)) \geqslant
(1/2)\cdot q(f_m(A), \partial D^{\,\prime})>\delta_1/2\,,$$
\begin{equation}\label{eq14}
q(f_m(|\gamma_{2, m}|))\geqslant q(z^{\,\prime}_m, f_m(y_0))
\geqslant (1/2)\cdot q(f_m(A),
\partial D^{\,\prime})>\delta_1/2\,,
\end{equation}
$$m\geqslant m_1>m_0\,.$$
Choose a chordal ball $U:=B_q(z_0, r_0),$ where $r_0> 0 $ and
$r_0<\delta_1/4,$ and $\delta_1$ is the number from
relations~(\ref{eq14}). Note that $f_m(|\gamma_{1, m}|)\cap
U\ne\varnothing\ne f_m(|\gamma_{1, m}|)\cap (D^{\,\prime}\setminus
U)$ for sufficiently large $m\in{\Bbb N},$ because
$q(f_m(|\gamma_{1, m}|))\geqslant \delta_1/2$ and $z_m\in
f_m(|\gamma_{1, m}|),$ $z_m\rightarrow z_0$ as $m\rightarrow\infty.$
Due to the same considerations $f_m(|\gamma_{2, m}|)\cap
U\ne\varnothing\ne f_m(|\gamma_{2, m}|)\cap (D^{\,\prime}\setminus
U).$ Since $f_m(|\gamma_{1, m}|)$ and $f_m(|\gamma_{2, m}|)$ are
continua, then by~\cite[Theorem~1.I.5.46]{Ku}
\begin{equation}\label{eq8AA}
f_m(|\gamma_{1, m}|)\cap \partial U\ne\varnothing, \quad
f_m(|\gamma_{2, m}|)\cap
\partial U\ne\varnothing\,.
\end{equation}
For a fixed $P>0,$ let $V\subset U$ be a neighborhood of the point
$z_0,$ corresponding to the definition of a weakly flat boundary,
that is, such that for any continua $E, F\subset D^{\,\prime}$ with
$E\cap
\partial U\ne\varnothing\ne E\cap \partial V$ and $F\cap \partial
U\ne\varnothing\ne F\cap \partial V$ the inequality
\begin{equation}\label{eq9AA}
M(\Gamma(E, F, D^{\,\prime}))>P
\end{equation}
holds. Note that for sufficiently large $m\in {\Bbb N}$
\begin{equation}\label{eq10AA}
f_m(|\gamma_{1, m}|)\cap \partial V\ne\varnothing, \quad
f_m(|\gamma_{2, m}|)\cap
\partial V\ne\varnothing\,.\end{equation}
Indeed, $z_m\in f_m(|\gamma_{1, m}|)$ and $z^{\,\prime}_m\in
f_m(|\gamma_{2, m}|),$ where $z_m, z^{\,\prime}_m\rightarrow z_0\in
V$ as $m\rightarrow\infty$ Therefore, $f_m(|\gamma_{1, m}|)\cap
V\ne\varnothing\ne f_m(|\gamma_{2, m}|)\cap V$ for large $m\in {\Bbb
N}.$ In addition, $q(V)\leqslant q(U)=2r_0<\delta_1/2$ and since
by~(\ref{eq14}) $q(f_m(|\gamma_{1, m}|))>\delta_1/2,$ then
$f_m(|\gamma_{1, m}|)\cap (D^{\,\prime}\setminus V)\ne\varnothing.$
Then $f_m(|\gamma_{1, m}|)\cap\partial V\ne\varnothing$
(see~\cite[Theorem~1.I.5.46]{Ku}). Similarly, $q(V)\leqslant
q(U)=2r_0<\delta_1/2$ and, since by~(\ref{eq14}) $q(f_m(|\gamma_{2,
m}|))>\delta,$ then $f_m(|\gamma_{2, m}|)\cap (D^{\,\prime}\setminus
V)\ne\varnothing.$ Now, by~\cite[Theorem~1.I.5.46]{Ku} we obtain
that $f_m(|\gamma_{1, m}|)\cap\partial V\ne\varnothing.$ Thus,
(\ref{eq10AA}) is proved.

\medskip
According to~(\ref{eq9AA}) and taking into account~(\ref{eq8AA}) and
(\ref{eq10AA}), we obtain that
$$M(f_m(\Gamma_m))=M(\Gamma(f_m(|\gamma_{1, m}|), f_m(|\gamma_{2, m}|),
D^{\,\prime}))>P\,,$$
which contradicts the inequality~(\ref{eq15}). The resulting
contradiction indicates that the original assumption made
in~(\ref{eq12}) is incorrect. The theorem is proved.~$\Box$

\section{Examples}

\begin{example}\label{ex1}
Let $D$ be the unit square from which the sequence of segments
$I_k=\{z=(x, y)\in {\Bbb R}^2: x=1/k,\,\,0<y<1/2\},$ $k=2,3,\ldots,$
is removed. See Figure~\ref{fig7} for this.
\begin{figure}
  \centering\includegraphics[width=450pt]{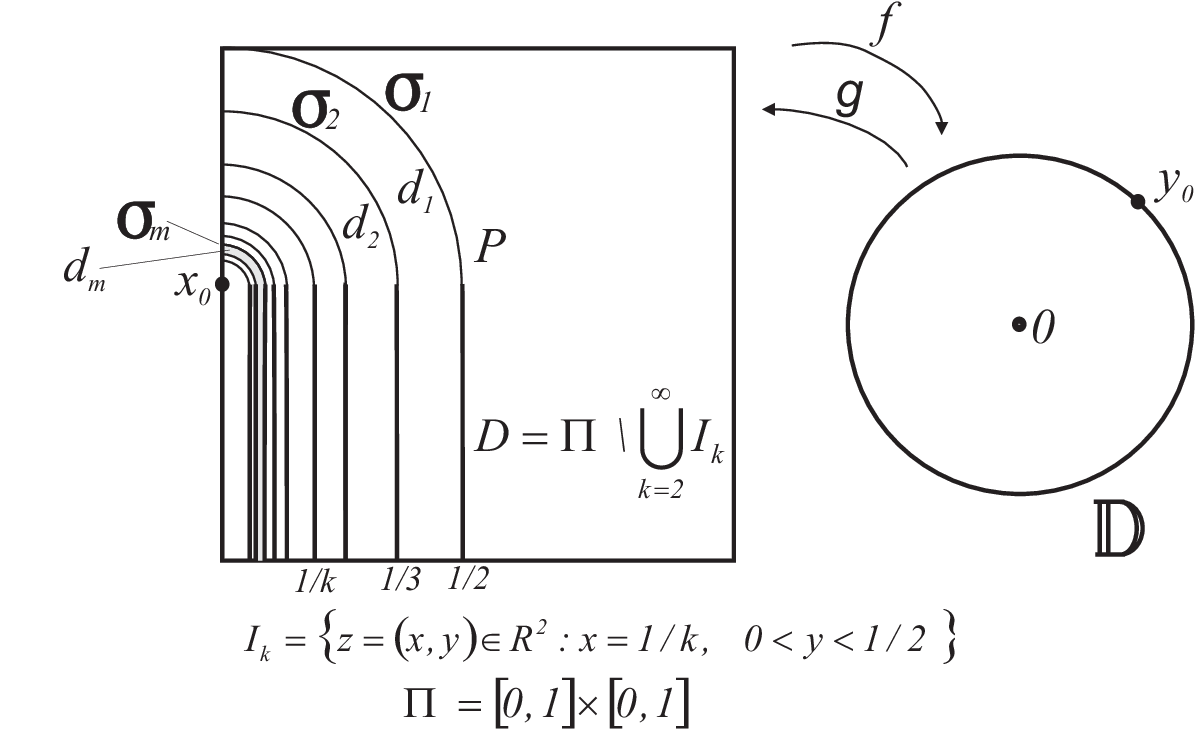}
  \caption{Illustration for Example~\ref{ex1}.}\label{fig7}
 \end{figure}
Consider the prime end $P$ in the domain $D,$ formed by cuts
$$\sigma_m=\left\{z=x_0+\frac{e^{i\varphi}}{m+1},\,\, x_0=(0, 1/2),\,\, 0\leqslant
\varphi\leqslant \pi/2\right\}, \quad m=1,2,\ldots, .$$
It can be shown that the end $P$ is really prime. According to the
Riemann mapping theorem, there exists a conformal mapping $g$ of the
unit disk ${\Bbb D}$ onto the domain $D$ and by the Caratheodory
theorem, a prime end $P$ corresponds to some point~$y_0\in\partial
{\Bbb D}$ so that $C(f, y_0)=I(P),$ see~\cite[Theorem~9.4]{CL}. It
follows that we may choose two sequences $z_k, w_k\in D,$
$k=1,2,\ldots ,$ such that $z_k, w_k\rightarrow P,$ $z_k\rightarrow
z_0$ and $w_k\rightarrow w_0$ as $k\rightarrow\infty,$ $z_0\ne w_0,$
while $f(z_k)\rightarrow y_0$ and $f(w_k)\rightarrow y_0$ as
$k\rightarrow\infty.$ Consequently, the mapping $f:=g^{\,-1}$ does
not have a continuous extension to the point $y_0$ in the pointwise
sense, although $g$ has a continuous extension
$\overline{g}:\overline{{\Bbb D}}\rightarrow \overline{D}_P.$

Consider another auxiliary family of mappings. As is known, linear
fractional automorphisms of the unit disk ${\Bbb D}\subset{\Bbb C}$
have the form
$$f(z)=e^{i\theta}\frac{z-a}{1-\overline{a}z}, \quad z\in {\Bbb
D},\quad a\in{\Bbb D},\quad\theta\in [0, 2\pi)\,.$$
We set, for example, $\theta=0,$ $a=1/n,$ $n=1,2,\ldots .$ In this
case, consider the family of mappings
$\widetilde{f}_n(z)=\frac{z-1/n}{1-z/n}=\frac{nz-1}{n-z}.$ Let
$\widetilde{A}=[0, 1/2].$ Then we obtain that
$\widetilde{f}_n(0)=-1/n\rightarrow 0$ and
$\widetilde{f}_n(1/2)=\frac{n-2}{2n-1}\rightarrow 1/2$ as
$n\rightarrow\infty.$ Thus, the mappings $\widetilde{f}_n$ satisfy
the condition $q(\widetilde{f}_n(\widetilde{A}))\geqslant\delta$
say, with $\delta=1/4.$

Now we put $f_n:=\widetilde{f}_n\circ f.$ Note that the mappings
$f_n$ are conformal; therefore, they satisfy the
estimate~(\ref{eq2*!}) for $Q\equiv 1$ at each point $x_0\in
\overline{D},$ see~~\cite[Theorem~1]{Pol}, cf.~\cite[Theorems~4.6
and 6.10]{MRSY$_1$} and \cite[Theorems~8.1 and 8.6]{MRSY$_3$}. Note
also that the mappings $f_n$ satisfy all the conditions of
Theorem~\ref{th2}, in particular, $q(f_n(A))\geqslant\delta$ with
$\delta=1/4,$ where $A:=f^{\,-1}(\widetilde{A}).$ Note that the unit
disk ${\Bbb D}$ has a weakly flat boundary by~\cite[Theorems~17.10
and~17.12]{Va$_1$}. By construction, the mappings
$g_n^{\,-1}:=f_n^{\,-1}$ do not even have a pointwise continuous
extension to $\partial {\Bbb D}$ in particular, the family of these
mappings is not equicontinuous as the family from $\overline{{\Bbb
D}}$ to $\overline{D}.$ Nevertheless, the extended family
$\overline{g}_n:\overline{\Bbb D}\rightarrow \overline{D}_P$ is
equicontinuous in terms of prime ends, as follows from
Theorem~\ref{th2}.

Put now $\theta=0,$ $a=(n-1)/n,$ $n=1,2,\ldots $ and
$\widetilde{f}^{-1}_n(z)=\frac{z-(n-1)/n}{1-z(n-1)/n}=\frac{nz-n+1}{n-nz+1}.$
Let $y_0=e^{i\theta},$ $0\leqslant \theta<2\pi.$ In this case, we
set $f_n:=\widetilde{f}_n \circ (e^{-i\theta}f).$ It is easy to
understand that the sequence $\widetilde{f}^{\,-1}_n$ is locally
uniformly converges to a constant function $-1$ in the unit disk. On
the other hand, we have the equality $\widetilde{f}^{\,-1}_n(1)=1,$
which immediately implies that the sequence $\widetilde{f}^{\,-1}_n$
is not equicontinuous at the point 1.

It follows that the sequence $f_n$ is also not equicontinuous at
point~$1.$ The reason for this is a violation of the
requirement~$q(\widetilde{f}_n(\widetilde{A}))\geqslant\delta.$
\end{example}

\begin{example}\label{ex2}
It is also easy to indicate a similar example of a family of
mappings with unbounded characteristic. Let $D$ be the domain
constructed in Example~\ref{ex2}. Then we put
$f_1(z)=\frac{1}{e\sqrt{2}}|z-1/2|.$ Note that $f_1$ maps $D$ onto a
domain $D_1$ lying in the ball $B(0, 1/e).$ Now we put
$f_2(z):=\frac{z}{|z|\log\frac{1}{|z|}}$ and $F(z)=(f_2\circ
f_1)(z).$ Using the technique outlined in the consideration
of~\cite[Proposition~6.3]{MRSY$_3$}, we may establish that $F$ is a
ring $Q$-homeomorphism in $\overline{D}$ with
$Q(z)=\log\frac{e\sqrt{2}}{|z-1/2|}.$ One can also prove that $Q\in
L^1(D).$ Note that $D_1$ is a simply connected domain, therefore, by
the Riemann theorem, it is possible to map it onto the unit disk
using some conformal mapping $f_3.$

Consider the family of mappings
$\widetilde{f}_n(z)=\frac{z-1/n}{1-z/n}=\frac{nz-1}{n-z}.$ We set
$F_n(z)=(\widetilde{f}_n\circ f_3\circ f_2\circ f_1)(z).$ Let
$\widetilde{A}=[0, 1/2].$ Then we obtain that
$\widetilde{f}_n(0)=-1/n\rightarrow 0$ and
$\widetilde{f}_n(1/2)=\frac{n-2}{2n-1}\rightarrow 1/2$ as
$n\rightarrow\infty.$ Thus, the mappings $\widetilde{f}_n$ satisfy
the condition $q(\widetilde{f}_n(\widetilde{A}))\geqslant\delta$
say, with $\delta=1/4.$ In this case, the mappings $F_n$ satisfy the
condition $q(F_n(A))\geqslant\delta$ say, with $\delta=1/4$ and
$A=(f_1^{\,-1}\circ f_2^{\,-1}\circ f_3^{\,-1})(\widetilde{A}).$

Since the modulus of families of paths does not change under
conformal transformations, the mappings $F_n$ are ring $Q$-maps in $
D,$ where $Q(z)=\log\frac{e\sqrt{2}}{|z-1/2|}$
(see~\cite[Theorem~8.1]{Va$_1$}). The mappings $G_n=F^{\,-1}_n$ do
not have a pointwise continuous extension to $\partial{\Bbb D},$
however, this extension is valid in the sense of prime ends. In
addition, the family of extended mappings
$\overline{G}_n:\overline{{\Bbb D}}\rightarrow \overline{D}_P,$
$n=1,2,\ldots ,$ is equicontinuous in $\overline{{\Bbb D}}$ by
Theorem~\ref{th2}.
\end{example}


\medskip
\medskip
{\bf \noindent Evgeny Sevost'yanov, Sergei Skvortsov, Nataliya Ilkevych} \\
Zhytomyr Ivan Franko State University,  \\
40 Bol'shaya Berdichevskaya Str., 10 008  Zhytomyr, UKRAINE \\
Phone: +38 -- (066) -- 959 50 34, \\
Email: esevostyanov2009@gmail.com

\end{document}